\documentclass{amsart}
\usepackage{amsmath,amssymb,enumerate,verbatim}

\usepackage{epsfig,fancyhdr,color}
\usepackage[applemac]{inputenc} 
  
\newcommand{\ds}{\displaystyle}
\newcommand{\tq}{\, \big| \, }

\DeclareMathOperator{\met}{Met}
\DeclareMathOperator{\aire}{Aire}
 \renewcommand{\r}{\mathbb{R}}

\newtheorem{theorem}{\rm\bf Théorème}[section]
\newtheorem{proposition}[theorem]{\rm\bf Proposition}
\newtheorem{lemme}[theorem]{\rm\bf Lemme}
\newtheorem{corollaire}[theorem]{\rm\bf Corollaire}

\theoremstyle{definition}

\theoremstyle{remark}

\def\interieur#1{\mathord{\mathop{\kern 0pt #1}\limits^\circ}}

\begin{document}
\title{Les surfaces à courbure intégrale bornée au sens d'Alexandrov}

\date{15 mai 2009}
\author{Marc Troyanov}
\address{\'Ecole Polytechnique F{\'e}d\'erale de
Lausanne, 1015 Lausanne, Bâtiment BCH - IGAT, Suisse}
\email{marc.troyanov@epfl.ch}


\keywords{surfaces d'Alexandrov,   courbure  bornée}


\begin{abstract}
Dans les années 1940-1970, Alexandrov et l'``\'Ecole de Leningrad" ont  
développé une théorie très riche des surfaces singulières. Il s'agit  
de surfaces topologiques, munie d'une métrique intrinsèque pour  
laquelle on peut définir une notion de courbure, qui est une mesure de  
Radon. Cette classe de surfaces a de bonnes propriétés de convergence  
et elle est remarquablement stable par rapport à diverses  
constructions géométriques  (recollements etc.). Elle englobe les  
surfaces polyédrales ainsi que les surfaces riemanniennes de classe  
$C^2$ ;  ces deux classes formant des parties denses de l'espace des  
surfaces d'Alexandrov. Toute surface singulière qu'on peut  
raisonnablement imaginer est une surface d'Alexandrov et de  
nombreuses propriétés géométriques des surfaces lisses s'étendent et  
se généralisent aux surfaces d'Alexandrov.

Le but de cet exposé est de donner une introduction non technique à  
la théorie d'Alexandrov, de donner des exemples et quelques-uns des  
faits fondamentaux de la théorie. Nous présenterons également un  
théorème de classification des surfaces (compactes) d'Alexandrov.

\medskip

\noindent AMS Mathematics Subject Classification:     \  53.c45
\end{abstract}

\maketitle

\section{Introduction}  

L'objet que nous appelons  ``surface''  n'est pas défini de façon univoque, et le choix d'un outillage plutôt qu'un autre nous incite a concevoir une surface comme un objet différentiable, polyédral, (semi)algébrique ou autre. Nous croyons définir un type d'objets et finissons par en étudier un autre. Souvenons-nous du célèbre épisode du mouchoir de Lebesgue. On raconte que, suite à un cours de Darboux où le maître démontrait que toute surface développable  (i.e. localement isométrique au plan)  dans l'espace euclidien $\mathbb{E}^3$ est une \emph{surface réglée} (c'est à dire que par chaque point de
cette surface  il passe un segment de droite de  $\mathbb{E}^3$ qui est contenu dans la surface), le jeune Lebesgue posa sur la table un mouchoir froissé : \emph{ce mouchoir est une surface développable; elle n'est pas réglée}. Darboux faisait de la géométrie différentielle, Lebesgue observe que dans la nature les surfaces sont rarement lisses.

\medskip

Pouvons-nous \emph{unifier}, et non \emph{opposer}, la géométrie différentielle lisse de Darboux et la géométrie froissée de Lebesgue ? Y a-t-il une théorie non triviale qui englobe toute les surfaces  que nous pouvons raisonnablement concevoir ?  Dans cet article, je me propose d'esquisser à traits grossiers, la solution que A.D. Alexandrov et ses élèves ont donnée à ce problème dans leur travaux sur les surfaces durant les années 1940-1970.

\medskip

Avant d'aller plus loin, voyons un premier problème qui semble facile à formuler. Nous appellerons \emph{surface lisse} une variété riemannienne $(S,g)$ de dimension $2$, et nous appellerons \emph{surface polyédrale}  un espace métrique qui est localement isométrique à un polyèdre (convexe ou non) de dimension $2$. Nous pouvons alors poser le

\medskip

\textbf{Problème 1. }
\begin{enumerate}[(a)]
  \item  Toute surface lisse est-elle limite d'une suite de surfaces polyédrales ?  
  \item  Toute surface polyédrale est-elle limite d'une suite de surfaces lisses ?
  \item  Quelles sont les invariants géométriques qui passent à la limite pour ces convergences ? 
\end{enumerate}

Le lecteur aura remarqué que la question n'est pas bien formulée : de quel type de convergence parle-t-on ? mais laissons cette question pour l'instant.

Le problème 1c) est d'une importance fondamentale, si nous voulons pouvoir passer d'un type de surface à un autre : du monde lisse au monde polyédral. Il faut comprendre quels invariants géométriques restent stables. Or ce problème est difficile et les premières observations donnent des réponses plutôt décourageantes.

\begin{itemize}
  \item[$\circ$] Le premier exemple est le \emph{lampion de H. A. Schwarz} \cite{Schwarz}. Il s'agit d'un exemple bien connu d'une suite de polyèdres (non convexes) $P_i$ qui converge vers un cylindre dans l'espace euclidien. L'aire $A$ de ce cylindre est finie mais l'aire de $P_i$ tend vers $+\infty$ avec $i$ (plus généralement on peut ajuster les paramètres des polyèdres pour obtenir une suite $\{P_i\}$ dont l'aire tend vers un nombre choisi arbitrairement
 dans l'intervalle $[A,+\infty]$, voir \cite[p. 117]{Gugg}). Ce phénomène et d'autres reliés est étudié dans l'article \cite{HPW}.

  \item[$\circ$] Un second exemple de non convergence est le groupe des isométries : une surface limite peut à priori avoir de nombreuses symétries lors même que les surfaces approchantes en sont dépourvues. Pour s'en rendre compte, on peut admirer la Géode de la Cité des Sciences à Paris. Cette belle construction est une approximation très convaincante d'une sphère par une surface polyédrale (triangulée). La plupart des sommets de ce polyèdre sont de valence six : il s'y rencontre 6 triangles. Mais la géode présente ici et là quelques ``défauts'' qui sont des sommets de valence 5. Ces défauts ont pour conséquences que la géode possède en réalité très peu de symétries alors qu'elle approxime la très symétrique sphère. Et il ne s'agit bien sûr pas d'une maladresse des ingénieurs qui ont conçu la géode, mais bien d'un problème intrinsèque à la géométrie (et en fait la topologie) de la sphère. 
\end{itemize}

\medskip

\textbf{Exercice :}  Décrire toutes les triangulations d'une sphère dont chaque sommet est de même valence.

\medskip

Ces exemples suffiront à nous convaincre que l'étude des convergences de surfaces nous confronte à d'intéressantes questions géométriques et qu'il peut-être utile de formaliser un peu mieux les questions du problème 1.

\medskip

\textbf{Problème 1'.}  Fixons une surface topologique $S$, on souhaite définir un espace $\mathcal{M}(S)$ contenant toutes les métriques ``raisonnables'' sur $S$ et munir cet espace d'une topologie pour laquelle les métriques polyédrales
sur $S$ et les métriques riemanniennes sur $S$ forment des parties denses. \\
On veut aussi décrire les invariants géométrique qui définissent des fonctions continues sur $\mathcal{M}(S)$.

\medskip

Les surfaces à courbure intégrale bornée au sens d'Alexandrov fournissent une réponse adéquate à ce problème.

\medskip

De nombreux résultats de la théorie des surfaces  à courbure intégrale bornée au sens d'Alexandrov ont été publiés en russe et n'ont jamais été traduits.
Les références les plus utiles en anglais sont les livres\footnote{Notons que le livre  \cite{AZ}
 d'Alexandrov et Zalgaller est disponible en ligne sur le site de l'American Mathematical Society.  
 Je recommande la lecture de l'introduction et du chapitre 1.} \cite{AZ,R93}. 
 Des travaux plus récents sur le sujets sont  \cite{BB1,BB2,BL,Shioya}. Je pense que le sujet devrait 
retrouver une nouvelle jeunesse notamment en raison des développements intenses de la géométrie différentielle discrète, voir 
 par exemple les livres \cite{Bobenko,GuYau}.


\section{Définition des surfaces d'Alexandrov}

Fixons une surface topologique $S$ que nous supposons pour simplifier compacte, sans bord et orientée\footnote{Rappelons qu'une surface est orientable si elle ne contient aucun ouvert homéomorphe au ruban de M\"obius.}. 

\medskip

\textbf{Définition} Une \emph{métrique à courbure intégrale bornée au sens d'Alexandrov} sur $S$ est la donnée d'une fonction 
continue 
$$
  d : S \times S  \to  \mathbb{R}
$$
telle que
\begin{enumerate}[(i)]
  \item $d$ est une distance et cette distance définit la topologie de variété sur $S$;
  \item l'espace métrique $(S,d)$ est \emph{convexe au sens de Menger}, c'est à dire qu'entre deux points $x,y\in S$   il existe toujours un point milieu $z$ tel que $d(x,z)=d(z,y) = \frac{1}{2}d(x,y)$;
  \item la distance $d$ sur $S$ est limite uniforme de distances associées à une suite de métriques riemanniennes sur $S$ pour laquelle l'intégrale de la valeur absolue de la courbure  est uniformément bornée.
\end{enumerate}
Il est bien connu que la condition (2) est équivalente à demander que la métrique $d$ est \emph{géodésique}, i.e. que toute paire de points $x,y \in S$ peut être reliée par une courbe de longueur $d(x,y)$ (voir \cite[théorème 2.6.2]{Papadopoulos2005}). La troisième condition justifie  l'appellation  ``courbure intégrale bornée''. Pour préciser ce qu'elle signifie, nous avons besoin de quelques rappels et préliminaires.

\medskip

Notons $\met(S)$ l'ensemble des métriques vérifiant les conditions (i) et (ii).  La distance uniforme entre deux éléments $d_1,d_2 \in \met(S)$ est définie par
\begin{equation}\label{ }
  D(d_1,d_2) = \sup \{ |d_1(x,y) - d_2(x,y)  | \, : \,  x,y \in S \}.
\end{equation}
On dit que la métrique $d\in \met(S)$ est \emph{riemannienne} si sont données une structure différentiable sur $S$ et une métrique riemannienne $g \in \Gamma(S^2T^*(S))$ telles que 
$$
 d(x,y) = \inf \int_0^1 \sqrt{g(\dot \gamma (t),\dot \gamma (t))} dt ,
$$
où l'infimum est pris sur l'ensemble de tous les chemins différentiables $\gamma : [0,1] \to S$ tels que $\gamma (0) = x$
et $\gamma (1) = y$.

\medskip

Rappelons maintenant ce que sont l'aire et la courbure de $g$. Sur tout ouvert simplement connexe  $U$ de $S$, on peut construire deux 
formes différentielles $\theta^1, \theta^2 \in \Omega^1(U)$ telles que $\theta^1 \wedge \theta^2$ est d'orientation positive  et
$$
 g = (\theta^1)^2 + (\theta^2)^2 = \theta^1\otimes \theta^1 + \theta^2\otimes \theta^2.
$$
Une paire de $1$-formes  $\theta^1, \theta^2$ vérifiant l'équation ci-dessus s'appelle un \emph{corepère mobile}.

\medskip

L'\emph{étoile de Hodge},  est l'endomorphisme $* : \Omega^k(S) \to  \Omega^{2-k}(S)$  ($k=0,1,2$)  vérifiant les conditions 
$$* (1) = \theta^1 \wedge \theta^2, \qquad  * (\theta^1 \wedge \theta^2) = 1, \qquad * \theta^1 =\theta^2, \qquad * \theta^2 = -\theta^1$$
pour tout corepère orthonormé $\theta^1,\theta^2$ d'orientation positive. 
La \emph{forme de connexion} associée au corepère $\theta^1, \theta^2$ est la $1$-forme  $\omega \in \Omega^1(U)$ définie par
\begin{equation}\label{eq.formconnection}
 \omega = -(*d\theta^1)\theta^1  -(*d\theta^2)\theta^2,
\end{equation}
c'est  l'unique $1$-forme sur $U$ vérifiant les équations de structure d'Elie Cartan :
$$
  \begin{cases}
  d\theta^1 &=  - \omega \wedge \theta^2 \\
  d\theta^2 &= \    \omega \wedge \theta^1.
  \end{cases}
$$

Un petit calcul nous donne le résultat suivant :

\begin{lemme}
Les formes différentielles $dA_g = \theta^1 \wedge \theta^2$ et $d\omega$ ne dépendent pas du choix du repère mobile $\theta^1, \theta^2 \in \Omega^1(U)$, elles sont donc globalement définies.
\end{lemme}

\medskip

On a  construit deux formes différentielles de degré $2$ sur $S$. Bien qu'elles ne soient pas exactes, nous les notons $d\omega$ et 
$$
 dA_g = \theta^1 \wedge \theta^2,
$$
nous y pensons comme deux mesures sur $S$ (on peut d'ailleurs montrer que  $dA_g$  est la mesure de Hausdorff  de $(S,d)$ en dimension 2).

\medskip

\textbf{Définition}  Nous appellerons $dA_g$ la \emph{mesure d'aire} sur $(S,g)$ et $d\omega$ la \emph{mesure de courbure}. 
La fonction $K : S \to \mathbb{R}$ définie par $K = *d\omega$ (i.e. $d\omega = K dA_g$) est la \emph{courbure Gaussienne} de $g$ (la courbure Gaussienne est donc la dérivée de Radon Nikodym de $d\omega$ par rapport à $dA_g$).
Remarquons  que la mesure de courbure $d\omega$  est invariante par homothétie, et rappelons qu'un résultat fondamental de la géométrie des surfaces,
la formule de Gauss-Bonnet, dit que l'intégrale de la courbure d'une surface close $(S,g)$ est égale à sa caractéristique d'Euler $\chi (S)$ :
$$
 \int_S d\omega = 2\pi \chi (S).
$$

Nous définissons aussi les mesures
$$
 d\omega^+ = K^+dA_g, \qquad d\omega^- = K^-dA_g \quad \text{et}  \quad  d |\omega| = |K|dA_g,
$$
où $K^+ = \max \{K,0\}$ et  $K^- = \max \{-K,0\}$. Observons que
$$
 d\omega = d\omega^+ -  d\omega^-,
  \qquad d |\omega| =d\omega^+ + d\omega^-.
$$

\medskip

Avec ces rappels, nous pouvons terminer la définition des métriques à courbure intégrale bornée. La condition (iii) signifie qu'il existe une suite de métriques riemanniennes $\{g_j\}$ sur $S$ telle que 
$$
 D(d_{g_j},d) \to 0  \qquad  \mathrm{et}   \qquad  \sup_{j \in \mathbb{N}} \int_S d\omega_{g_j}^+ < \infty.
$$

La formule de Gauss-Bonnet entraîne que $\int_S d|\omega_j|$ est aussi bornée et le théorème de Banach-Alaoglu nous dit alors qu'il existe une mesure de Radon $d\omega$ sur $S$ et une sous-suite  $\{g_{j'}\}$ de  $\{g_j\}$ telle que
$$
  d\omega_{g_{j'}}  \to d\omega
$$
faiblement\footnote{Dire que $d\omega_{g_{j'}}  \to d\omega$ faiblement signifie que $\int_S f d\omega_{g_{j'}}  \to \int_S f d\omega$ pour toute fonction continue $f\in C(S)$. Cette condition entraîne que $\omega_{g_{j'}}(E)  \to \omega(E)$ pour tout ensemble mesurable $E\subset S$ dont la frontière
$\mathrm{Fr}(E)$ vérifie $\omega(\mathrm{Fr}(E))=0$, voir \cite[p. 134]{Doob}.}.

\begin{theorem}
 La mesure de Radon $d\omega = \lim_{j'} d\omega_{g_{j'}}$ est bien définie sur la surface métrique $(S,d)$, elle ne dépend pas de  la  suite $\{g_j\}$ choisie.
\end{theorem}

Ce résultat est un cas particulier du théorème \ref{th.conv.ALex} ci-dessous. Il n'est pas banal puisqu'à priori la courbure dépend des dérivées jusqu'à l'ordre $2$ du tenseur métrique et que la convergence considérée est une convergence uniforme.

\medskip

\textbf{Définition}  La mesure limite $d\omega $ s'appelle la \emph{mesure de courbure} de la surface d'Alexandrov $(S,d)$.

\medskip

Observons que la formule de Gauss-Bonnet est vérifiée pour toute surface d'Alexandrov à cause de la continuité des intégrales pour la convergence faible des mesures. Le corollaire \ref{cor.realdmu} ci-dessous nous dit qu'il n'y a pas d'autre condition : \emph{toute mesure de Radon sur une surface close compatible avec la formule de Gauss-Bonnet est une mesure de courbure pour une métrique d'Alexandrov sur} $S$.

\medskip

\textbf{Remarque} Ce qui est impliqué par ce théorème est que l'on peut construire la mesure $d\omega$ directement à partir de la métrique $d$ sur la surface $S$ (sans passer par des approximations riemanniennes). L'idée de la construction est de trianguler la surface par des triangles géodésiques et de compter dans chaque triangle l'excès angulaire\footnote{de façon naïve, l'\emph{excès angulaire}  d'un triangle est  égal à (somme des angles  - $\pi$); mais la notion d'angle n'est pas vraiment bien définie et Alexandrov la remplace part une notion d'\emph{angle supérieur} (\emph{upperangle}).}, puis de faire tendre la maille de la triangulation vers $0$. La construction rappelle celle de Carathéodory pour produire une mesure à partir d'une mesure extérieure,  voir   \cite[chap. 5]{AZ}.

\medskip

La définition des surfaces  à courbure intégrale bornée au sens d'Alexandrov que nous avons donnée n'est pas la définition donnée par Alexandrov dans \cite{AZ}, et c'est donc un théorème.
La manière dont Alexandrov définit les surfaces  à courbure intégrale bornée est synthétique, i.e. purement métrique. Il demande que l'excès angulaire total de toute famille de triangles 
simples dans $S$ dont les intérieurs sont deux à deux disjoints soit borné uniformément. Cette définition est expliquée en détails aux pages 4--6 de \cite{AZ}, voir aussi \cite[71-74]{R93}.

\section{Surfaces d'Alexandrov obtenues par recollement}

Il nous faut des exemples. Considérons une famille finie $\{ T_1, T_2, \dots , T_m\}$ de triangles Riemanniens, c'est à dire de variétés homéomorphes à un disques, munies de métriques riemanniennes dont le bord est de classe $C^2$ par morceaux avec trois points anguleux. 

Recollons ces triangles selon une triangulation prescrite de notre surface $S$ en recollant les paires de côtés par des isométries. On obtient un espace métrique (de longueur) homéomorphe à $S$ qui est une surface d'Alexandrov dont la mesure de courbure est donnée par 
$$
 d\omega = d\omega_0 + d\omega _1 + d\omega_2,
$$
où $d\omega_2$ est absolument continue par rapport à la mesure d'aire et donné par
$$
 d\omega_2 = K dA
$$
à l'intérieur de chaque triangle. La mesure $d\omega_0$ est une mesure discrète supportée par les sommets de la triangulation, telle que pour chaque sommet $p$, on a
$$
 \omega_0(\{p\}) = 2\pi - (\text{la somme des angles des triangles $T_i$ incidents au sommet $p$} ).
$$
Finalement la mesure $d\omega _1$ est supportée par les arêtes de la triangulation, et sur chaque arête $a = T_i\cap T_j$ , on a
$$
 d\omega_1 = (k^+ - k^-) ds
$$
où $k^+$ et $k^-$ sont les courbures géodésiques de $a$ vues dans chacun des deux triangles adjacents à $a$
orientés de façon cohérente avec un choix d'orientation de l'arête.

\medskip

Une \emph{surface polyédrale} est une surface obtenue par recollement de triangles euclidiens (au sens classique, i.e. dont le bord est formé de trois segments
de droites). La mesure de courbure est alors concentrée aux sommets de la triangulation et au voisinage de chaque sommet, la surface est localement isométrique à un cône euclidien. Pour cette raison une surface polyédrale s'appelle aussi une \emph{surface euclidienne à singularités coniques}. Inversément, toute surface d'Alexandrov dont la mesure de courbure est discrète est une surface polyédrale.

\medskip

Voyons quelques exemples. Considérons la surface d'un cube (que nous pouvons trianguler si nous le désirons). Les faces sont plates et donc $d\omega_2 = 0$, les arêtes sont géodésiques vues de chaque face incidentes et  donc $d\omega_1 = 0$. Les huit sommets sont incidents chacun à 3 angles de $\frac{\pi}{2}$, nous avons donc
$$
 \omega (p) = 2\pi - 3\frac{\pi}{2} = \frac{\pi}{2},
$$
en chacun des huit sommets. 
On vérifie directement la formule de Gauss-Bonnet :
$$
 \int_S d\omega = \int_S d\omega_0 = 8\times \frac{\pi}{2} = 4 \pi = 2\pi \chi(S^2).
$$
Comme second exemple, considérons une boîte de conserve. C'est un cylindre euclidien de rayon $r$ et de hauteur $h$ complété par son fond et son couvercle qui sont deux disques euclidien $D_1,D_2$ de rayon $r$. On peut trianguler topologiquement les différentes pièces de cette boîte de conserve si on le désire. Le cylindre et les deux disques sont des surfaces plates, donc la courbure se concentre sur les deux cercles bordant les deux disques. Ces cercles sont géodésiques comme courbes dans le cylindre, et comme courbes dans les disques $D_i$ ils sont de courbure géodésique constante $k = \frac{1}{r}$.
On a donc $d\omega_0= d\omega_2 = 0$ et 
$$
  d\omega_1 = \frac{1}{r} \left. ds \right|_{\partial D_1} + \frac{1}{r} \left. ds \right|_{\partial D_2}.
$$
On vérifie à nouveau la formule de Gauss-Bonnet :
$$
 \int_S d\omega = \int_S d\omega_1 = \frac{1}{r} \text{Long} (\partial D_1) +  \frac{1}{r} \text{Long} (\partial D_2)
 =  \frac{1}{r} (2\pi r + 2\pi r) = 4 \pi.
$$

\section{Structure conforme et uniformisation des surfaces lisses}

Dans ce paragraphe, nous considérons des métriques riemanniennes lisses. Rappelons qu'une métrique $\tilde{g}$ sur la surface différentiable $S$ est une \emph{déformation conforme} de $g$ s'il existe une fonction $u : S \to \r$ telle que $\tilde{g} = e^{2u}g$. Si $u$ est une constante, on dit que la métrique $\tilde{g}$ est \emph{homothétique} à $g$.

Un sytème de coordonnées $(x,y)$ sur un ouvert $U \subset S$ est dit \emph{conforme} pour la métrique $g$  s'il existe une fonction $\rho : U \to \r$ telle que 
$$
 g = \rho(x,y) (dx^2+dy^2).
$$
Ces coordonnées s'appellent aussi des\emph{ coordonnées isothermes}. 

\medskip

\textbf{Définition. } Le \emph{laplacien} de $u$ par rapport à la métrique $g$ est défini par 
$$\Delta_g u = -* d*du,$$
 c'est à dire $\Delta_g u dA = -d*du$.
Le laplacien est un opérateur elliptique,  il s'écrit en coordonnées
$$
 \Delta_g u = -\frac{1}{\sqrt{\det (g_{ij})} } \sum_{\mu,\nu = 1}^2 \frac{\partial}{\partial x_{\mu}}\left(g^{\mu\nu} {\sqrt{\det (g_{ij})} } \cdot \frac{\partial u}{\partial x_{\nu}}\right).
$$
Dans le cas des coordonnées conformes, si $g = \rho(x,y) (dx^2+dy^2)$, alors 
$$
 \Delta_g u = -\frac{1}{\rho(x,y)}\left(\frac{\partial^2 u}{\partial x^2} + \frac{\partial^2 u}{\partial y^2}\right).
$$

\medskip

\begin{proposition}
 Toute surface riemannienne lisse admet des coordonnées conformes au voisinage de chacun de ses point.
\end{proposition} 

\textbf{Preuve.}  Soit $U \subset S$ un ouvert dans lequel est défini un corepère mobile $\theta^1,\theta^2\in \Omega^1(U)$ tel que $g = (\theta^1)^2 + (\theta^1)^2$. Alors $\tilde{\theta}^1 = e^u \theta^1$, $\tilde{\theta}^2 = e^u \theta^2$ est un repère mobile pour la métrique $\tilde{g} = e^{2u}g$, et un examen des équations de structures montre facilement que les formes de connexions $\omega$ et $\widetilde{\omega}$ sont reliées par l'équation
$$
  \widetilde{\omega} = \omega - *du.
$$
Par conséquent 
\begin{equation}\label{eq.chgt.courbure0}
 d\widetilde{\omega} = d\omega - d*du.
\end{equation}

Les mesures d'aires des métriques $g$ et $\tilde{g}$ sont données par $dA = \theta^1 \wedge \theta^2$ et $d\widetilde{A} = \widetilde{\theta}^1 \wedge \widetilde{\theta}^2 = e^{2u}dA$.
L'équation (\ref{eq.chgt.courbure0})  peut donc s'écrire
\begin{equation}\label{eq.chgt.courbure}
   \tilde{K}e^{2u} = K + \Delta_g u.
\end{equation}

Choisissons une fonction $u : U \to \r$ telle que $\Delta_g u = -K$, alors  $\tilde{g} = e^{2u}g$ est une métrique plate (i.e. de courbure nulle) sur $U$ et on peut
donc trouver au voisinage de chaque point de $U$ des coordonnées euclidiennes telles que  $\tilde{g} = dx^2+dy^2$. Dans ce voisinage on a donc
$$
  g = e^{-2u}(dx^2+dy^2).
$$
\qed

\bigskip

\textbf{Remarques}  \  \textbf{1)}  Ce théorème s'appelle classiquement le ``théorème d'existence de coordonnées isothermes''. Il a une histoire intéressante, il remonte à Gauss dans le cas où la métrique $g$ est analytique. En 1914-16, Korn et Lichtenchtein \cite{K,L} démontrent ce résultats en supposant que $g$ est seulement H\"older continue, et Chern \cite{Chern}  a simplifié leur preuve en 1955. Les travaux d'Ahlfors et Bers ont permis d'étendre ce résultat au cas des surfaces dont la métrique est seulement mesurable, sous certaines conditions. 

\smallskip

\textbf{2)}  L'équation (\ref{eq.chgt.courbure}) joue un rôle fondamental en géométrie des surfaces. Remarquons qu'en intégrant cette équation, on trouve
$$
 \int_S  \tilde{K}d\tilde{A} = \int_S  \tilde{K}e^{2u}dA = \int_S  KdA +  \int_S  \Delta_g udA,
$$
or   $\int_S  \Delta_g udA = -\int_S d * du = 0$ par la formule de Stokes, donc
$$
 \int_S  \tilde{K}d\tilde{A} =   \int_S  KdA.
$$
Cette identité est  compatible avec la formule de Gauss-Bonnet.

\begin{corollaire} \label{cor1.strcomp}
 Toute métrique riemannienne $g$ sur une surface différentiable orientée définit une structure complexe.
\end{corollaire}

\textbf{Preuve.} On peut construire un atlas orienté sur $(S,g)$ dont tous les sytèmes de coordonnées sont conformes. Les changements
de coordonnées vérifient alors les conditions de Cauchy-Riemann et sont donc holomorphes.

\qed

\medskip

{\small On peut aussi aussi observer que l'étoile de Hodge définit par dualité une structure presque complexe $J : TS \to TS$, or en dimension réelle 2, toute 
structure presque complexe est intégrable. Cet argument prouve à la fois la proposition et le corollaire.}

\medskip

Nous pouvons maintenant citer le théorème d'Uniformisation de Poincaré-Koebe.

\begin{theorem}
 Toute métrique riemannienne $g$ sur une surface compacte est une déformation conforme d'une métrique $h$ à courbure constante.
\end{theorem}

\textbf{Une esquisse de preuve.}  Soit $(S,g)$ une surface riemannienne compacte et sans bord. Supposons d'abord que $\chi (S) = 0$, alors  $\int_S  KdA =0$.
On sait alors qu'il existe une solution $u\in C^{\infty}(S)$ de l'équation 
\begin{equation}
 \Delta_g u = -K.
\end{equation}
L'équation (\ref{eq.chgt.courbure}), entraîne que $h = e^{2u}g$ est une métrique de courbure identiquement nulle.

Dans le cas où  $\chi (S) < 0$, le raisonnement est semblable. On résoud l'équation non linéaire
\begin{equation}\label{eq.Knonlin}
    \Delta_g u -  \tilde{K}e^{2u} = 1.
\end{equation}
La métrique $h = e^{2u}g$ est alors une métrique de courbure constante $-1$; l'équation (\ref{eq.Knonlin}) peut se résoudre par une méthode variationnelle,
c'est ce qui a été fait en 1969  par Melvyn Berger \cite{B, Aubin}.  On peut aussi résoudre cette équation par des méthodes de flot de Ricci.
Dans le cas où $\chi (S) > 0$, l'argument ne marche pas pour des raisons techniques assez subtiles. On peut contourner la difficulté en introduisant une singularité,
voir \cite[p. 624]{Troyanov1990}.

\qed  

\medskip

Dans la suite de cet article, nous normaliserons les métriques à courbure constante en demandant que $K=-1$ si $\chi(S) < 0$ et 
$K= +1$ si $\chi(S) > 0$. Si  $\chi(S) = 0$ alors $K=0$  et on normalise la métrique en demandant que $(S,h)$  soit d'aire $=1$.

\section{Fonction de Green et potentiel}

L'inverse du laplacien est donné par la fonction de Green (voir \cite{Aubin, deRham}) :
\begin{theorem}
Soit $(S,h)$ une surface riemannienne lisse, compacte et sans bord. Alors il existe une unique fonction 
$G : S\times S \to  \r \cup \{ +\infty\}$ vérifiant les conditions suivantes~:
\begin{enumerate}[(a)]
  \item $G$ est $C^{\infty}$ sur $ S\times S \setminus \{(x,x) \tq x \in S \}$;
  \item $G(x,y) = G(y,x)$;
  \item $|G(x,y)| \leq C \cdot (1+ |\log d(x,y)|)$;
  \item $\ds \int_S G(x,y) dA_h (y) = 0$;
  \item  Pour toute fonction   $u \in C^2(S)$, on a 
  $$
   u(x) =  \int_S G(x,y) \Delta u (y)  dA_h (y) + \frac{1}{\aire (S)} \int_S  u (y)  dA_h (y) .
  $$
\end{enumerate}
\end{theorem}

\begin{proposition}
 Soit $\mu$ une mesure de Radon d'intégrale nulle sur $S$, alors la fonction
\begin{equation}\label{eq.pot}
  u(x) =  \int_S G(x,y) d\mu (y)
\end{equation}
 vérifie l'équation $\Delta u = \mu$ au sens faible.
\end{proposition}

\textbf{Preuve} Soit $\varphi$ une fonction test, i.e. $\varphi \in C^{\infty}(M)$, et notons 
$\overline{\varphi } =  \frac{1}{\aire (S)} \int_S   \varphi (y)   dA_h (y)$. Alors on a
\begin{align*}
 \int_S  \varphi (x) \cdot  \Delta u(x) \;  dA(x) & =  \int_S  \Delta \varphi (x) \cdot u(x) \; dA(x) 
  \\  & =   \int_S \int_S  \Delta \varphi (x)  G(x,y)  \; d\mu (y)dA(x) 
  \\  & =   \int_S \int_S  G(x,y) \Delta \varphi (x)   \; dA(x)  d\mu (y)
  \\  & =    \int_S \left( \varphi (y)  - \overline{\varphi }\right)  d\mu (y)
  \\  & =    \int_S  \varphi (y)  d\mu (y).
\end{align*}
\qed

\medskip

\textbf{Définition} On dit que la fonction $u$ définie par l'équation (\ref{eq.pot}) est le \emph{potentiel} de la mesure $d\mu$ relativement à la métrique $h$. 

\medskip

Nous pouvons dire un certain nombre de choses d'une fonction $u$ qui est potentiel d'une mesure. Tout d'abord $u$
est d'intégrale nulle, et $u$ est différence de deux fonctions sous-harmoniques. La régularité de $u$ peut se décrire
par 
$$
 u \in \bigcap_{p<2} W^{1,p}(S),
$$
(voir \cite[théorème 9.1]{Stampacchia}), de plus on a
$$
 \sup_{0< \varepsilon \leq 1} \sqrt{\epsilon} \|  \nabla u \|_{L^{2-\varepsilon}(S)} \leq C(h, |\mu| (S))
$$
(voir \cite[théorème 2]{Iw}).
La fonction $u$ est en outre approximatevement différentiable presque partout et on a 
$$
  \| \text{ap} \nabla u \|^*_{L^{2,\infty}(S)} \leq C(h, |\mu| (S))
$$
où $L^{2,\infty}(S)$ est un espace de Lorentz (voir \cite[théorème 2]{DHM}).

 \medskip

Dans la suite, on notera $V(S,h)$ l'ensemble des fonctions $u$ sur $S$  telles que $\mu = \Delta_h u$ est une mesure.
Pour toute fonction $u\in V(S,h)$ et tous $x,y \in S$, on note
\begin{equation} \label{def.qdist}
 d_{h,u}(x,y) = \inf  \left\{ \int_0^1 e^{u(\alpha(t))} |\dot \alpha(t)|_hdt  \tq  \alpha \in \mathcal{C}(x,y) \right\}
\end{equation}
où $\mathcal{C}(x,y)$ est l'ensemble des  chemins $\alpha : [0,1] \to S$, de classe $C^1$ tels que $\alpha(0) = x$ et $\alpha (1) = y$. Il est clair que $d_{h,u}$ est une pseudo-métrique, i.e. que $d_{h,u}$ est symétrique et vérifie l'inégalité du triangle, de plus $0 \leq d_{h,u}(x,y) \leq \infty$ pour tous $x,y \in S$. 

\begin{proposition}
 La pseudo-métrique $d_{h,u}$ est séparante, i.e. $d_{h,u}(x,y) > 0$ si $x\neq y$. De plus $d_{h,u}(x,y) < \infty$ pour toute paire de points $x,y\in S$ telle que $\mu (\{ x\}) < 2\pi$ et $\mu(\{ y\}) < 2\pi$ où $\mu$ est la mesure $\Delta_h u$.
\end{proposition}

Ce résultat a été observé par Youri Reshetnyak. Il se déduit facilement des lemmes 4.1 et 4.2 de \cite{Troyanov1991}.
Un point $x$ tel que $\Delta_h u(\{ x\}) = 2\pi$ s'appelle un \emph{cusp}. Il peut se situer à distance finie ou infinie.

\medskip

Résumons ce que nous disent les considérations précédentes : si $(S,h)$ est une surface riemannienne compacte et $\mu$ est une mesure de Radon sur $S$ d'intégrale nulle et telle que $\mu (\{ x\}) < 2\pi$ pour tout $x\in S$, et $u$ est le potentiel de $\mu$, alors 
$g = e^{2u}h$ est une métrique riemannienne singulière sur $S$ pour laquelle la pseudo-distance associée (définie par (\ref{def.qdist})) est une vraie distance. 

On verra bientôt que $(S, d_{h,u})$ est une surface à courbure intégrale bornée au sens d'Alexandrov.

\section{Théorèmes de convergence}

Un résultat fondamental dit que la mesure de courbure d'une surface à courbure intégrale bornée au sens d'Alexandrov dépend continûment de la métrique, c'est le théorème 6 page 240 de \cite[chap. 7]{AZ}.

\begin{theorem}\label{th.conv.ALex}
 Soit $(S,d)$ une  surface à courbure intégrale bornée au sens d'Alexandrov  compacte et $\{ d_j\}$ une suite de métrique 
 à courbure intégrale bornée convergeant vers $d$ dans la topologie uniforme. Alors la mesure de courbure 
 de $(S,d)$ est limite faible des mesures de courbure de $\{(S,d_j)\}$, i.e. 
 $$
  D(d_j,d) \to 0 \ \Longrightarrow  \ \int_S f(x) d\omega_j (x) \to  \int_S f(x) d\omega (x)
 $$
 pour toute fonction continue $f$ sur $S$, où $d\omega_j$ est la mesure de courbure de $(S,d_j)$ et 
 $d\omega$ est la mesure de courbure de $(S,d)$.
\end{theorem}

En 1960, Youri Reshetnyak démontre à l'inverse que la métrique dépend continûment de la mesure de courbure, 
à condition qu'il n'y ait aucun cusp et que la structure conforme soit fixée, c'est le Théorème III de \cite{R60a} :

\begin{theorem}\label{th.conv.Resh}
 Soit $(S,h)$ une surface riemannienne lisse et $\{d\mu^+_n\}$, $\{d\mu^-_n\}$ deux suites de mesures de
 Radon sur $S$ convergeant faiblement vers les mesures $d\mu^+ = \lim_{n\to \infty} d\mu^+_n$ et
 $d\mu^- = \lim_{n\to \infty} d\mu^-_n$. Supposons que $\mu^+ (\{ p\}) < 2\pi$ pour tout point $p\in S$.
 Notons $u_n$ le potentiel de $d\mu_n = d\mu^+_n + d\mu^-_n$ et $u$ le potentiel de
 $d\mu = d\mu^+ + d\mu^-$, alors 
 $$
   d_{h,u_n} \to d_{h,u}
 $$
 dans la topologie uniforme.
\end{theorem}

\medskip

\begin{corollaire}
Soit $(S,h)$ une surface riemannienne compacte lisse et $\mu$ une mesure de Radon  sur $S$ 
d'intégrale nulle et telle que $\mu (\{ x\}) < 2\pi$ pour tout $x\in S$. Soit $u$ le potentiel de $\mu$, alors la métrique 
$d_{h,u}$ définie par (\ref{def.qdist}) est une métrique à courbure intégrale bornée au sens d'Alexandrov.
La mesure de courbure de $(S,d_{h,u})$ est donnée par
\begin{equation}\label{eq.chgt.courbure2}
 d\omega = K_hdA_h + d\mu.
\end{equation}
\end{corollaire}
 
\emph{Preuve}  Donnons-nous une suite de mesures lisses $d\mu_j = \varphi_j dA_h$ qui converge faiblement vers $d\mu$. Soit $u_j$ le potentiel de $d\mu_j$, alors $g_j = e^{2u_j}h$ est une métrique lisse
et la distance associée $d_j$ converge vers $d_{h,u}$ pour la topologie uniforme par le théorème \ref{th.conv.Resh}. Donc $(S,d_{h,u})$ est une surface   à courbure intégrale bornée au sens d'Alexandrov par définition.
 L'identité (\ref{eq.chgt.courbure2}) se déduit alors de (\ref{eq.chgt.courbure}) par passage à la limite.
 
 \qed
 
 \medskip
 
 Le résultat suivant dit que la structure conforme est déterminée par la métrique $d_{h,u}$.
 
  \medskip
  
\begin{theorem}\label{th,rigiditeconforme}
 Soient $(S,h)$ et $(S',h')$ deux surfaces riemanniennes compactes lisses et $u\in V(S,h)$, $u'\in V(S',h')$.
 Notons  $d_{h,u}$, respectivement $d_{h',u'}$, les métriques d'Alexandrov associées à $e^uh$ et
 $e^{u'}h'$. \\
 Si $f : (S,d_{h,u}) \to (S',d_{h',u'})$ est une isométrie, alors $f$ est une transformation conforme de $(S,h)$
 vers $(S',h')$.
 \end{theorem}
 
La preuve peut se déduire d'un théorème de Menchoff (1937)  en théorie des applications quasi-conformes qui dit que tout homéomorphisme $1$-quasi-conforme entre deux surfaces est une application conforme (voir \cite{M,G}).

\section{Structure conforme et uniformisation des surfaces d'Alexandrov}

Reshetnyak \cite{R60b} (voir aussi \cite{Huber}) a  démontré à partir du théorème \ref{th.conv.Resh} que toute métrique d'Alexandrov sur une surface $S$ détermine une structure conforme.

\begin{theorem}\label{th.Exresh}
 Soit $(S,d)$ une surface à courbure intégrale bornée au sens d'Alexandrov compacte sans cusp. Alors 
 il existe une métrique riemannienne $h$ et une fonction $u\in V(S,h)$ telle que 
 $$
   d = d_{h,u}.
 $$
\end{theorem}

\begin{corollaire}
 Toute métrique  à courbure intégrale bornée $d$ au sens d'Alexandrov sans cusp sur une surface compacte $S$ orientée définit une structure complexe
 sur cette surface.
\end{corollaire}

\textbf{Preuve.}  Le théorème précédent nous dit que $d = d_{h,u}$ où $h$ est  une métrique riemannienne et $u\in V(S,h)$. Le corollaire 
\ref{cor1.strcomp} permet de définir une structure complexe sur $(S,h)$ et le théorème \ref{th,rigiditeconforme} entraîne que cette  structure 
complexe  est uniquement déterminée.

\qed

\medskip

Voyons quelques exemples de surfaces à courbure intégrale bornée en représentation conforme :

\begin{enumerate}
  \item Soit $V$ un c\^one euclidien d'angle $\theta$, alors $V$
est isom\'etrique \`a $\mathbb{C}$ muni de la m\'etrique
$$ds^2 = |z|^{2\beta}|dz|^2 \, ,$$
o\`u $\displaystyle \beta = (\frac{\theta}{2\pi} -1)$ (cf. \cite[prop.1]{Troyanov1986}).
 Cette m\'etrique est
donc de classe $L^p$ pour $1 < p < -1/\beta $ si $\beta < 0$ (i.e. si
$\theta < 2\pi$) et de classe $L^{\infty}$ si $\beta > 0$.
Sa courbure est la mesure
$$d\omega = -2\pi \beta \cdot \delta _0 \, ,$$
o\`u $\delta _0 $ est la mesure de Dirac en $0$.
  \item Consid\'erons la surface $S$ obtenue en recollant une
demi-sph\`ere de rayon 1 \`a un demi-cylindre de rayon 1 le long de leur
bord par une isom\'etrie.
Alors $S$ est isom\'etrique \`a $(\mathbb{C},ds^2)$ o\`u $ds^2 = \rho (z) |dz|^2$
avec
$$\rho (z) = \begin{cases}
\frac{4}{(1+|z|^2)^2} \, , & \text{si $|z| \le 1$ ;} \\ \ds
\frac{1}{|z|^2} \  , & \text{si $|z| \ge 1 $ .} 
\end{cases} 
$$
On voit donc que cette surface est de classe $C^{1,1}$. Sa mesure de courbure est absolument continue : 
$d\omega = KdA$ avec
$$
K (z) = \begin{cases}
1 \,  & \text{si $|z| < 1$ } ;\\
0 \,   & \text{si $|z| > 1 $ .} 
\end{cases} 
$$
  \item Soit $S$ la surface obtenue en recollant 
deux disques euclidiens de rayons $1$ par isom\'etrie le long de leur bord.
Alors $S$ est topologiquement une sph\`ere munie d'une m\'etrique plate sur
le compl\'ementaire d'un cercle $\Sigma$ form\'e de points singuliers.
On montre facilement que $S$ est isom\'etrique \`a $\mathbb{C} \cup \{ \infty
\}$ muni de la m\'etrique $\rho (z) |dz|^2$ o\`u
$$
\rho (z) =  \begin{cases}
1 \,  & \text{si $|z| \le 1$ ;} \\ \ds
|z|^{-4} \,   & \text{si $|z| \ge 1 $ .} 
\end{cases} 
 $$
Cette m\'etrique est donc de classe $C^{0,1}$ (Lipschitz). Sa mesure de courbure 
est donn\'ee par 
$$d\omega = 2 \cdot ds _{| \Sigma} , $$
o\`u $ds _{| \Sigma}$ est la mesure donn\'ee par la longueur riemannienne le
long de $\Sigma$.
\item Soit $P$ la pseudo-sph\`ere de Beltrami (\`a laquelle on
ajoute son point \`a l'infini). Alors $P$ est isom\'etrique au disque $D =
\{ z : |z| < 1\}$ muni de la m\'etrique
$$ds^2 = \frac{|dz|^2}{(|z|\log (z))^2} \, .$$
Cette m\'etrique est de classe $L^1$, l'origine est un cusp à distance infinie et la mesure de courbure est 
$$d\omega = -dA + 2\pi \cdot \delta _0 .$$
\end{enumerate}

{\small Dans  ces exemples, la mesure de courbure se calcule
en consid\'erant une suite de m\'etriques lisses $\{ g_i \}$
approximant convenablement la m\'etrique singuli\`ere $ds^2$. On peut en
effet v\'erifier dans chaque exemple que $K_jdA_j$ converge faiblement vers
$d\omega$. Il est facile de vérifier la formule de Gauss-Bonnet dans chacun de ces exemples}.

\bigskip

Résumons : 
Une métrique d'Alexandrov sans cusp sur une surface compacte $S$  détermine les données suivantes
 \begin{enumerate}[(i)]
  \item Une structure conforme sur $S$;
  \item La mesure de courbure $d\omega$.
\end{enumerate}

Inversement :
\begin{theorem}\label{th.clsf}
Pour toute structure conforme sur $S$ et toute mesure de Radon  $d\omega$ telle que  
 $\int_S d\omega = 2\pi \chi (S)$ et $\omega (\{x\}) < 2\pi$ pour tout $x$,  il existe une métrique d'Alexandrov sur $S$ associée à cette structure conforme et dont la mesure de courbure est égale à $d\omega$.
Cette métrique d'Alexandrov est unique à une homothétie près.
\end{theorem}

\textbf{Preuve} Donnons-nous une structure conforme sur $S$, que nous représentons par une métrique riemannienne $h$ à courbure constante. Soit $d\mu =  d\omega - K_hdA_h$ et $u\in V(S,h)$ le potentiel de 
$d\mu$. Alors la métrique $d_{h,u}$ a les propriétés voulues.

Pour prouver l'unicité, on considère une autre métrique d'Alexandrov $d'$ sur $S$. Par le théorème de Reshetnyak \ref{th.Exresh}, on sait qu'il existe une  métrique riemannienne $h'$ sur $S$ et une 
fonction $u'\in V(S',h')$ telles que $d' = d_{h',u'}$.
Le théorème  \ref{th,rigiditeconforme} entraîne que $h$ et $h'$ sont conformément équivalente, i.e. il existe une fonction $v \in C^{\infty}(S)$ telle que $h'  = e^{2v}h$. Quitte à remplacer $u'$ par $u'+v$, on peut donc supposer que $h=h'$. On a donc $d' = d_{h,u'}$ dont la mesure de courbure est $d\omega$, par conséquent
$$
 \Delta_h u' = d\omega  - K_hdA_h =  \Delta_h u.
$$
Ainsi $\Delta_h (u'- u) = 0$ et $(u'-u)$ est donc constante.

\qed 

\medskip

Ce résultat peut-être vu comme un théorème de classification des surfaces d'Alexandrov. 
Soit $S$ une surface compacte orientée, notons $\mathcal{M}_0(S)$ l'espace des métriques à courbure intégrale bornée sur $S$ sans cusp, $\mathcal{C}(S)$ l'espace des structures conformes
et $\mathcal{R}_{2\pi}(S)$ l'espace des mesures de Radon  $d\omega$ sur $S$ telles
que
$$
 \int_S d\omega = 2\pi\chi(S) \qquad \text{et} \qquad \omega(\{x\}) < 2\pi
 \ \text{ pour tout } x \in S.
$$
Alors le théorème précédent dit que
$$
 \mathcal{M}_0(S) =  \mathcal{C}(S) \times \mathcal{R}_{2\pi}(S)\times \r_+.
$$
(où terme $\r_+$ contrôlant le facteur d'homothétie).
Notons que l'espace des mesures de Radon sur $S$ est lui-même un espace métrisable localement complet (cf. par exemple \cite{Doob}) et que l'espace des structures complexes est 
bien compris, notamment via le point de vue de la théorie de Teichmüller.

\medskip

Voyons quelques conséquences de ce théorème. La première conséquence est le résultat suivant sur les surfaces euclidiennes à singularités coniques (voir \cite{Troyanov1986}).

\begin{corollaire}
Soit $S$ une surface close, $x_1, x_2,\cdots, x_n$ des points de $S$ et $\theta_1, \theta_2,\cdots \theta_n > 0$. Supposons que $\sum_i(2\pi-\theta_i)=2\pi\chi(S)$, alors pour toute structure conforme sur $S$, il existe une métrique polyédrale conforme sur $S$ ayant en $x_i$ une singularité  conique  d'angle 
$\theta_i$ $(i=1, \cdots,n)$. Cette métrique est unique à homothétie près.
\end{corollaire}

\textbf{Preuve} Il suffit d'appliquer le théorème précédent à la mesure discrète
$$
  d\omega = \sum_i(2\pi-\theta_i) \delta_{x_i}.
$$
\qed

\medskip 

Plus généralement, on a

\begin{corollaire}\label{cor.realdmu}
 Toute   mesure de Radon  $d\omega$ sur $S$ telle que  
 $\int_S d\omega = 2\pi \chi (S)$ et $\omega (\{x\}) < 2\pi$ pour tout $x$  est la mesure de courbure
 d'une métrique d'Alexandrov.
\end{corollaire}

\medskip

\begin{corollaire}
 Toute surface d'Alexandrov est limite d'une suite de surfaces polyédrales.
\end{corollaire}

\textbf{Preuve} Il suffit d'approximer la mesure de courbure de la  surface d'Alexandrov donnée par une
suite de mesures discrètes.

\qed

Rappelons que selon notre définition, toute surface d'Alexandrov est trivialement limite d'une suite de surfaces riemaniennes.


\end{document}